# An Autoregressive Model with Semi-stable Marginals

**S Satheesh**

NEELOLPALAM, S. N. Park Road
Trichur – 680 004, **India.**
*ssatheesh1963@yahoo.co.in*

**E Sandhya**

Department of Statistics, Prajyoti Niketan College
Pudukkad, Trichur – 680 301, **India.**
*esandhya@hotmail.com*

**Abstract.** The family of semi-stable laws is shown to be semi-selfdecomposable. Thus they qualify to model stationary first order autoregressive schemes. A connection between these autoregressive schemes with semi-stable marginals and semi-selfsimilar processes is given.

**1. Introduction.**

The notion of semi-stable laws is by Levy and for various aspects of it see, Satheesh and Sandhya (2006*a*). Recently Maejima and Naito (1998) introduced the notion of semi-selfdecomposable (SSD) laws and Satheesh and Sandhya (2005) observed that SSD laws can generate marginally stationary additive first order autoregressive (AR(1)) schemes. Here we prove that semi-stable($a,b$) laws are SSD($b$). Now we are interested in describing an AR(1) model with semi-stable marginals. We need the following notions.

**Definition.1.1** (Pillai, 1971). A characteristic function (CF) $f$ is semi-stable($a,b$) if $\forall \; u \in \mathbf{R}$ and for some $0<|b|<1<a$, $f(u) = \{f(bu)\}^a$. Here $a$ and $b$ are connected by $a|b|^{\alpha} = 1$, $\alpha \in (0,2]$.

**Definition.1.2** (Maejima and Sato, 1999). A CF $f$ that is infinitely divisible (ID) is semi-stable($a,b$) if $f(u) = \{f(bu)\}^a$, $\forall \; u \in \mathbf{R}$, for some $a \in (0,1) \cup (1,\infty)$ and $b>0$. Here also $a$ and $b$ are connected by $ab^{\alpha} = 1$, $\alpha \in (0,2]$.

**Definition.1.3** (Maejima and Naito, 1998). A probability distribution $\mu$ on $\mathbf{R}$, with CF $f$ is SSD($b$) if for some $b \in (0,1)$ there exists a CF $f_o$ that is ID and

$$f(u) = f(bu) f_o(u), \forall\ u \in \mathbf{R}.$$

**Definition.1.4** (Maejima and Sato, 1999). A process $\{X(t), t \geq 0\}$ is semi-selfsimilar if for some $a>0$ there is a unique $H>0$ such that $\{X(at)\} \stackrel{d}{=} \{a^H X(t)\}$. We write $\{X(t)\}$ is $(a,H)$-semi-selfsimilar where $a$ is called the epoch and $H$ the exponent of the semi-selfsimilar process. If this relation holds for any $a>0$, then $\{X(t)\}$ is H-selfsimilar.

A process $\{X(t), t \geq 0\}$, $X(o) = 0$, having stationary and independent increments is a Levy process. A Levy process $X(t)$ such that $X(1)$ is semi-stable will be called a semi-stable process. Since semi-stable laws are ID (we will discuss this) these Levy processes are well defined. Clearly, selfsimilar processes are processes that are invariant in distribution under suitable scaling of time and space. From Maejima and Sato (1999) we also have: A Levy process $\{X(t)\}$ is semi-selfsimilar (selfsimilar) *iff* the distribution of $X(1)$ is semi-stable (stable). Here the notions of semi-stability and semi-selfsimilarity are considered in the strict sense only.

The additive AR(1) scheme that we consider here is described by the sequence of *r.v*s $\{X_n\}$, if there exists an innovation sequence $\{\varepsilon_n\}$ of *i.i.d r.v*s satisfying

$$X_n = bX_{n-1} + \varepsilon_n, \forall\ n>0 \text{ integer and some } 0<b<1. \tag{1}$$

Satheesh and Sandhya (2005) have also discussed methods to construct SSD laws, its implication in subordination of Levy processes and defined integer-valued SSD($b$) laws and corresponding integer-valued AR(1) model.

We compare the two definitions of semi-stable laws in section.2 and show that semi-stable($a,b$) laws are SSD($b$). Stationary AR(1) schemes with semi-stable marginals are then discussed connecting it to semi-selfsimilar processes.

**2. Results.**

**Remark.2.1** Apart from the range of the parameter $b$, definitions 1.1 and 1.2 differ on an important aspect. Definition.1.2 assumes the CF $f$ to be ID which is not there in definition.1.1. In general a CF $f(u)$, $u \in \mathbf{R}$ is complex-valued and may have a zero point. But no way is known how to define $\{f(u)\}^a$ for $u$ beyond a zero point except when $a>0$ is an integer. Possibly, this is why Maejima and Sato (1999) assume "*f* is ID" and further state (remark.4.1) that the assumption "*f* is ID" is not needed if the parameter *a*>0 is an integer. To circumvent this situation in definition.1.1 of Pillai (1971) we need assume that *f* has no zeroes. Otherwise this definition in terms of CFs, is meaningless. However, it is worth noting that Pillai (1971) showed that every semi-stable law is partially attracted to itself. Hence every semi-stable law must be ID. Also, such a situation will not be there in defining semi-stable laws on **R**$^+$ using Laplace transforms or if we describe it as the weak limit of appropriately normed partial sums.

**Remark.2.2** Since *a* and *b* are connected by $ab^\alpha = 1$, $\alpha \in (0,2]$ if *a*>1 then *b*<1 and if *a*<1 then *b*>1. Further, without loss of generality we may restrict the range of *a* or *b* as done in remark.2.3. Hence setting *a*>1 we have *b*<1. Now the description in Maejima and Sato (1999) gives a subclass of the semi-stable laws defined in Pillai (1971).

**Remark.2.3** Without loss of generality we may consider the range of *a* as 0<*a*<1 in the description of semi-selfsimilarity because it is equivalent to $\{(a^{-1})^H X(t)\} \stackrel{d}{=} \{X(a^{-1}t)\}$ and thus the whole range of *a*>0 is covered.

In the following discussion we will use the parameter range 0<*b*<1 for semi-stable(*a*,*b*) laws since we are discussing its connection to SSD(*b*) laws and subsequently to AR(1) model.

**Theorem.2.1** Semi-stable(*a*,*b*), 0<*b*<1, family of laws is SSD(*b*).

**Proof.** Let *f* be the CF of a semi-stable(*a*,*b*) law. By virtue of remark.2.1 *f* is also ID and hence;

$f(u) = \{f(bu)\}^a$, $\forall u \in \mathbf{R}$ and some 0<*b*<1<*a*.

  $= f(bu)\{f(bu)\}^{a-1}$.

Here the second factor is also ID and hence by definition.1.3 *f* is SSD(*b*).

Thus by Satheesh and Sandhya (2005) the semi-stable(*a*,*b*) family of laws qualify to model marginally stationary additive AR(1) schemes. Rather than just prescribing $\varepsilon_n$ to have the CF $\{f(bu)\}^{a-1}$, let us take a different look at this. Here we make use of the semi-selfsimilar processes we have briefly discussed.

**Theorem.2.2** (A corollary to the theorem.4.1 in Maejima and Sato (1999)). $\{X(t)\}$ is semi-stable$(a,b)$ Levy *iff* $\{X(t)\}$ is $(b^{-\alpha}, \frac{1}{\alpha})$-semi-selfsimilar.

**Proof.** Since $X(t)$ is semi-stable$(a,b)$ Levy, $bX(at) \stackrel{d}{=} X(t)$ or $X(at) \stackrel{d}{=} \frac{1}{b}X(t)$. Since $ab^\alpha = 1$, $\alpha \in (0,2]$, $\frac{1}{b} = a^{\frac{1}{\alpha}}$ and hence $X(at) \stackrel{d}{=} a^{\frac{1}{\alpha}} X(t)$. Thus $X(t)$ is $(a, \frac{1}{\alpha})$-semi-selfsimilar or, $X(t)$ is $(b^{-\alpha}, \frac{1}{\alpha})$-semi-selfsimilar. Converse easily follows.

**Theorem.2.3** Let $\{Z(t), t \geq 0\}$ be a Levy process, $X_0 \stackrel{d}{=} Z(1)$ and $\varepsilon_n \stackrel{d}{=} bZ(b^{-\alpha}-1)$, $\forall n$ in (1). Then (1) is marginally stationary with semi-stable$(a,b)$ marginals if $\{Z(t)\}$ is $(b^{-\alpha}, \frac{1}{\alpha})$-semi-selfsimilar. Conversely, $\{Z(t)\}$ is $(b^{-\alpha}, \frac{1}{\alpha})$-semi-selfsimilar and the marginals are semi-stable$(a,b)$ if (1) is marginally stationary.

**Proof.** Notice that if the CF of $Z(1)$ is $f(u)$ then that of $bZ(b^{-\alpha}-1)$ is $\{f(bu)\}^{b^{-\alpha}-1}$. Further if $\{Z(t)\}$ is $(b^{-\alpha}, \frac{1}{\alpha})$-semi-selfsimilar then $\{f(\frac{u}{b})\} = \{f(u)\}^{b^{-\alpha}}$ and so $\{f(bu)\}^{b^{-\alpha}} = f(u)$. Now under the given assumptions at $n=1$,

$$f_1(u) = f(bu) \{f(bu)\}^{b^{-\alpha}-1} = \{f(bu)\}^{b^{-\alpha}} = f(u).$$

Thus on iteration (1) is marginally stationary with semi-stable$(a,b)$ marginals.

Conversely, let (1) is marginally stationary. Then at $n=1$,

$$f(u) = f(bu) \{f(bu)\}^{b^{-\alpha}-1} = \{f(bu)\}^{b^{-\alpha}}.$$

Hence the marginals and $Z(1)$ are semi-stable$(a,b)$ and so $\{Z(t)\}$ is $(b^{-\alpha}, \frac{1}{\alpha})$-semi-selfsimilar. Thus the proof is complete.

**Concluding Remarks.** As corollaries to theorems 2.2 and 2.3 we have: $\{X(t)\}$ is stable Levy *iff* $\{X(t)\}$ is $\frac{1}{\alpha}$-selfsimilar. Similarly, (1) is marginally stationary with stable marginals if $\{Z(t)\}$ is selfsimilar. Conversely, $\{Z(t)\}$ is selfsimilar and the marginals are stable if (1) is marginally stationary. Notice that if the condition describing semi-stability and semi-selfsimilarity is true for all $a>0$ (more precisely, for two reals $a_1$ and $a_2$ such that $ln(a_1)/ ln(a_2)$ is irrational) then the condition describes stability and selfsimilarity. Satheesh and Sandhya (2006*b*) has considered the integer-valued analogue of the results in this paper.

**Acknowledgement.** Authors thank the referee for pointing out a mistake in the earlier version and for the comments that lead to remark.2.1.